\documentclass{article}

\usepackage[dvips]{graphicx}
\usepackage{amssymb}
\usepackage{amsmath}

\usepackage{caption}
\usepackage{hyperref}
\usepackage{arcs}
\usepackage {xcolor} 
\usepackage{subfigure}
\usepackage{authblk}
\usepackage{calrsfs}

\newtheorem{theorem}{Theorem}

\usepackage{graphicx} % Required for inserting images

\newenvironment{AMS}{\small\bf 2020 AMS subject classification: }{} 

\title{Nonsingularity of unsymmetric Kansa matrices: random collocation by MultiQuadrics\\ and Inverse MultiQuadrics}
\author[1]{R. Cavoretto}
\affil[1,3]{University of Torino, Italy}
\author[2]{F. Dell'Accio}
\affil[2]{University of Calabria, Rende (CS), Italy}
\author[3]{A. De Rossi}
\author[4]{A. Sommariva} 
\affil[4,5]{University of Padova, Italy}
\author[5]{M. Vianello}
\date{\today}

\begin{document}

\maketitle

\begin{abstract}

Unisolvence of unsymmetric Kansa collocation is still a substantially open problem. We prove that Kansa matrices with MultiQuadrics and Inverse MultiQuadrics for the Dirichlet problem of the Poisson equation are almost surely nonsingular, when the collocation points are chosen by any continuous random distribution in the domain interior and arbitrarily on its boundary.  
\end{abstract}

\vskip0.5cm
\noindent
\begin{AMS}
{\rm 65D12, 65N35.}
\end{AMS}
\vskip1.5cm

\section{Introduction}

Unsymmetric Kansa collocation has become over the years one of the most adopted meshless methods by RBF for the numerical solution of PDEs in a variety of engineering and scientific problems; cf., e.g., \cite{K86,K90} and \cite{CDR20,CFC14,LOS06,SW06} with the references therein. On the other hand, in the popular textbook \cite{F07} one reads: ``{\em Since the numerical experiments by Hon and Schaback show that Kansa's method cannot be well-posed for arbitrary center locations, it is now an open question to find sufficient conditions on the center locations that guarantee invertibility of the Kansa matrix}''. Indeed, Hon and Schaback \cite{HS01} proved that there are ``rare'' point configurations that make Kansa matrices singular, but so far unisolvence of unsymmetric Kansa collocation has remained a substantially open problem. 

In a quite recent paper \cite{DASV24}, unisolvence of {\em random} Kansa collocation by Thin-Plate Splines (TPS) has been proved for the Poisson equation, on 2-dimensional domains with analytic boundary. The key properties in the proof are radiality of the Laplacian and the fact that the TPS basis functions are analytic up to their center, the latter used also in other recent papers on interpolation unisolvence of TPS without polynomial addition, cf. \cite{BSV24,DASV23}. Though this result represents a first step towards unisolvence, it has a number of restrictions, besides the fact that the differential operator is the pure Laplacian: the RBF kind (indeed, the most usual approach to Kansa method is with MultiQuadrics), the dimension, the  boundary regularity. 

In the present paper, still resorting to the key property of analiticity of the basis, but this time with the presence of complex singularities, we prove that unsymmetric Kansa matrices with MultiQuadrics (MQ) and Inverse MultiQuadrics (IMQ) for the Dirichlet problem of the Poisson equation are almost surely invertible (in any dimension), when the collocation points are chosen by any continuous random distribution in the domain interior and arbitrarily on its boundary. 

\section{Unisolvence of random MQ and IMQ Kansa collocation}

In this paper we study unisolvence of Kansa collocation by (scaled) MultiQuadrics (MQ)
\begin{equation}\label{MQ}
\phi(r)=\phi_\varepsilon(r)=\sqrt{1+(\varepsilon r)^2}\;,
\end{equation}
and Inverse MultiQuadrics (IMQ)
\begin{equation}\label{IMQ}
\phi(r)=\phi_\varepsilon(r)=\frac{1}{\sqrt{1+(\varepsilon r)^2}}\;,
\end{equation}
which are both analytic in $\mathbb{R}$. The scale $\varepsilon >0$ represents the so-called \emph{shape parameter} associated with RBF \cite{F07,LS23}. We consider the Poisson equation with Dirichlet boundary conditions (cf. e.g. \cite{E98})
\begin{equation}\label{poisson_eq}
\left\{
\begin{array}{l}
\Delta u(P)=f(P)\;,\;P\in \Omega\;,
\\
u(P)=g(P)\;,\;P\in \partial \Omega\;,
\end{array} 
\right .
\end{equation}
where $\Omega\subset \mathbb{R}^d$ is a bounded domain (connected open set), $P=(x_1,\dots,x_d)$ and $\Delta=\partial ^2/\partial x_1^2+\dots+\partial ^2/\partial x_d^2$ is the Laplacian. Differently from \cite{DASV24}, where we studied Kansa discretization by TPS in $\mathbb{R}^2$ and we assumed that the boundary is a curve possessing an analytic parametrization, here we do not make any restrictive assumption on $\partial\Omega$, except for the usual ones that guarantee well-posedness and regularity of the solution (like e.g. that the boundary is Lipschitz, cf. e.g. \cite{S98} with the references therein). The main reason is that for the discretization of the boundary conditions, with MQ we can resort to a classical result by Micchelli \cite{M86} on interpolation unisolvence by any set of distinct points, result achieved also with IMQ since they are strictly positive definite. 

Unsymmetric Kansa collocation (see e.g. \cite{CDR20b,F07,HS01,K86,LOS06,SW06,W05}) consists in seeking a function
\begin{equation} \label{u_N}
u_N(P)=\sum_{j=1}^n{c_j\,\phi_j(P)}+\sum_{k=1}^m{d_k\,\psi_k(P)}\;,\;\;N=n+m\;,
\end{equation}
where
\begin{equation} \label{phij}
\phi_j(P)=\phi(\|P-P_j\|_2)\;,\;\;\{P_1,\dots,P_n\}\subset \Omega\;,
\end{equation}
\begin{equation} \label{psik}
\psi_k(P)=\phi(\|P-Q_k\|_2)\;,\;\;\{Q_1,\dots,Q_m\}\subset 
\partial\Omega\;,
\end{equation}
such that
\begin{equation}\label{poisson_disc}
\left\{
\begin{array}{l}
\Delta u_N(P_i)=f(P_i)\;,\;i=1,\ldots,n
\\
u_N(Q_h)=g(Q_h)\;,\;h=1,\ldots,m\;.
\end{array} 
\right .
\end{equation}

The following properties will be used below. Defining $\phi_A(P)=\phi(\|P-A\|)$, we have $\phi_A(A)=1$ and  $\phi_A(B)=\phi_B(A)$. Moreover 
$\Delta\phi_A(B)=\Delta\phi_B(A)$ and $\Delta\phi_A(A)=2\varepsilon^2$. Indeed, the Laplacian in polar coordinates 
centered at $A$ (cf. e.g. \cite[Ch.2]{E98}, \cite[App.D]{F07}) is the radial function
\begin{equation} \label{lapl}
\Delta\phi_A={\frac{\partial^2 \phi}{\partial^2 r}}+{\frac{1}{r}} \frac{\partial \phi}{\partial r}
\end{equation}
and thus for $\phi(r)=(1+(\varepsilon r)^2)^s,\, s\in {\mathbb R}\setminus\{0\}$, we get 
\begin{equation} \label{Lp}
\Delta\phi_A=4\varepsilon^2s\,(1+(\varepsilon r)^2)^{s-2}(1+s(\varepsilon r)^2)\;\;,
\end{equation}
that is for MQ ($s=1/2$)
\begin{equation} \label{LIMQ}
\Delta\phi_A=\varepsilon^2\,\frac{2+(\varepsilon r)^2}{(1+(\varepsilon r)^2)^{3/2}}\;\;,
\end{equation}
and for IMQ ($s=-1/2$)
\begin{equation} \label{LIMQ}
\Delta\phi_A=-\varepsilon^2\,\frac{2-(\varepsilon r)^2}{(1+(\varepsilon r)^2)^{5/2}}\;\;.
\end{equation}

Kansa collocation can be rewritten in matrix form as
\begin{equation} \label{Kansa-system}
\left(\begin{array}{cc}
\Delta\Phi & \Delta\Psi\\
\\
\Phi & \Psi
\end{array} \right) 
\left(\begin{array}{c}
\mathbf{c}\\
\\
\mathbf{d}
\end{array} \right) 
=
\left(\begin{array}{c}
\mathbf{f}\\
\\
\mathbf{g}
\end{array} \right) 
\end{equation}
where the block matrix is 
$$
K_n=K_{n,m}(\{P_i\},\{Q_h\})=\left(\begin{array}{cc}
\Delta\Phi & \Delta\Psi\\
\\
\Phi & \Psi
\end{array} \right) 
$$
$$
=\left(\begin{array} {ccccccc}
c & \cdots &  \cdots & \Delta\phi_n(P_1)
& \Delta\psi_1(P_1) & \cdots & \Delta\psi_m(P_1) \\
\\
\vdots & \ddots  &  & \vdots & \vdots  & \cdots & \vdots\\

\vdots &  & \ddots & \vdots & \vdots  & \cdots & \vdots\\
\\
\Delta\phi_1(P_n) & \cdots & \cdots & c
& \Delta\psi_1(P_n) & \cdots & \Delta\psi_m(P_n) \\
\\
\phi_1(Q_1) & \cdots & \cdots & \phi_n(Q_1)
& 1 & \cdots & \psi_m(Q_1) \\
\\
\vdots & \cdots & \cdots & \vdots & \vdots  & \ddots & \vdots\\
\\
\phi_1(Q_m) & \cdots & \cdots & \phi_n(Q_m)
& \psi_1(Q_m) & \cdots & 1\\
\\
\end{array} \right)
$$
with $c=2\varepsilon^2$ for MQ and $c=-2\varepsilon^2$ for IMQ, 
 and $\textbf{f}=\{f(P_i)\}_{i=1,\ldots,n}$, ${\textbf{g}}=\{g(Q_h)\}_{h=1,\ldots,m}$.
We are now ready to state and prove our main result. 

\vskip0.2cm
\begin{theorem}
Let $K_n$ be the MQ or IMQ Kansa collocation matrix defined above, where $\{Q_h\}$ is any fixed set of $m$ distinct points on $\partial\Omega$, and  $\{P_i\}$ is a sequence of i.i.d. (independent and identically distributed) random points in $\Omega$ with respect to any probability density $\sigma \in L^1_+(\Omega)$. 
Then for every $m \geq 1$ and for every $n\geq 0$ the matrix $K_n$ is a.s. (almost surely) nonsingular.
\end{theorem}
\vskip0.1cm

Before proving the theorem, we recall that the construction of i.i.d. random sequences with respect to any probability density can be accomplished (via the uniform distribution) by the well-known acceptance-rejection method, cf. e.g. \cite{Flury90}. Though uniform random points could be the most natural choice for collocation, the possibility of adopting other distributions could be interesting whenever it is known that the solution has steep gradients, or other regions where it is useful to increase the discretization density. 
\vskip0.2cm 
\noindent
{\bf Proof of Theorem 1.}
The proof proceeds by (complete) induction on $n$. For $n=0$ the collocation matrix coincides with the $m\times m$ interpolation matrix on the boundary discretization points, which is (deterministically) nonsingular. For IMQ this is a consequence of their positive definiteness (cf. e.g. \cite{F07,W05}), while for MQ  
this comes from a classical result by Micchelli on conditionally positive definite RBF of order 1, cf. \cite{M86}. 
For the inductive step, we define the augmented matrix 
$$
\tilde{K}(P)=\left(\begin{array} {cccccccc}
c & \cdots &  \cdots & \Delta\phi_n(P_1) & \Delta\phi_1(P)
& \Delta\psi_1(P_1) & \cdots & \Delta\psi_m(P_1) \\
\\
\vdots & \ddots  &  & \vdots & \vdots  & 
\vdots & \cdots & \vdots \\
\vdots &  & \ddots & \vdots & \vdots  & \vdots & \cdots & \vdots\\
\\
\Delta\phi_1(P_n) & \cdots & \cdots & c
& \Delta\phi_n(P) & \Delta\psi_1(P_n) & \cdots & \Delta\psi_m(P_n) \\
\\
\Delta\phi_1(P) & \cdots & 
\cdots & \Delta\phi_n(P) & c
& \Delta\psi_1(P) & \cdots & \Delta\psi_m(P) \\
\\
\phi_1(Q_1) & \cdots & \cdots & \phi_n(Q_1)
& \psi_1(P) & 1 &  \cdots & \psi_m(Q_1) \\
\\
\vdots & \cdots & \cdots & \vdots & \vdots  & \vdots & \ddots & \vdots \\
\\
\phi_1(Q_m) & \cdots & \cdots & \phi_n(Q_m)
& \psi_m(P) & \psi_1(Q_m) & \cdots & 1 \\
\\

\end{array} \right)
$$

Observe that in this case $K_{n+1}=\tilde{K}(P_{n+1})$
since $\psi_k(P_{n+1})=\phi_{n+1}(Q_k)$ and 
$\Delta\phi_j(P_i)=\Delta\phi_i(P_j)$.

To compute the determinant, developing $\det(V(P))$ by Laplace's rule on the $(n+1)$-row we have
\begin{equation} \label{delta}
\delta(P)=\mbox{det}(\tilde{K}(P))= -\mbox{det}(K_{n-1})(\Delta\phi_n(P))^2+\alpha(P)\Delta\phi_n(P)+\beta(P)
\end{equation}
where 
\begin{equation} \label{alpha-beta}
\alpha\in \mbox{span}\{\Delta\phi_j,\psi_k,\Delta\psi_k\,;\,1\leq j\leq n-1\,,\,1\leq k\leq m\}\;,
\end{equation}
$$
\beta\in \mbox{span}\{\Delta\phi_i\Delta\phi_j,\Delta\phi_i\Delta\psi_h,\psi_k\Delta\phi_i,\psi_k\Delta\psi_h\,;\,1\leq i,j\leq n-1\,,\,1\leq k,h\leq m\}\;.
$$

Notice that $\delta$ is a real analytic function in $\Omega$, because such are all the functions involved in its definition by linear combinations and products, and real analytic functions form a function algebra \cite{KP02}. 
We claim that $\delta(P)$ is almost surely not identically zero in $\Omega$. 
Indeed, if $\delta$ were identically zero in $\Omega$, it would be identically zero also in $\mathbb{R}^2$, since the zero set of a not identically zero real analytic function must have null Lebesgue measure (cf. \cite{M20} for an elementary proof) whereas $meas(\Omega)>0$. Then taking the line $P(t)=P_n+tv$ where $v=(v_1,\dots,v_d)$ is a given unit vector, we obtain that 
the real univariate function $\delta(P(t))$ would be identically zero for $t\in \mathbb{R}$. Consequently, its analytic extension  to the complex plane, say $\delta(P(z))$, would also be  identically zero for $z\in \mathbb{C}$.
Observe now that by (\ref{Lp})-(\ref{LIMQ})
$$
\Delta\phi_n(P(z))=4\varepsilon^2s\,(1+(\varepsilon z)^2)^{s-2}(1+s(\varepsilon z)^2)\;,\;s=1/2\;\;\mbox{or}\;\;s=-1/2\;,
$$ 
has two branching points in $z=\pm i/\varepsilon$ (we take in (\ref{MQ})-(\ref{IMQ}) the branch of the square root that is positive on the positive reals), and $(\Delta\phi_n(P(z)))^2$ has a pole there, of order 3 for MQ and of order 5 for IMQ. 
On the other hand, the functions $\alpha(P(z))$ and $\beta(P(z))$ are analytic at $z=\pm i/\varepsilon$. In fact, if $A$ is one of the collocation points different from $P_n$, that is $A\in \{Q_h\}\subset \partial \Omega$ or $A\in \{P_k, \,k\neq n\}\subset \Omega$, 
we first observe that $\|P(z)-A\|_2^2=\|P_n+zv-A\|_2^2$ has to be seen as the complex extension of the corresponding real function, hence not the complex 2-norm but the sum of the squares of the complex components. 
Then the complex numbers $$
1+\varepsilon^2\|P(\pm i/\varepsilon)-A\|_2^2=1+\varepsilon^2\sum_{j=1}^d(P_n\pm iv/\varepsilon -A)_j^2
$$
$$
=1+\varepsilon^2\sum_{j=1}^d[(P_n -A)_j^2\pm 2i(P-A)_jv_j/\varepsilon-v_j^2/\varepsilon^2]
$$
(recalling that $\|v\|_2=1$) 
have a.s. positive real part, 
namely $\|P_n-A\|_2^2>0$, since $P_n$ is a.s. distinct from $A$, and thus the complex functions corresponding to the chosen branch of the complex square root 
$$
\left(1+\varepsilon^2\sum_{j=1}^d(P_n+zv-A)_j^2\right)^{\pm 1/2}
$$
are both analytic at $z=\pm i/\varepsilon$.
This means that $\phi_A(P(z))$ and $\Delta\phi_A(P(z))$ are analytic at $z=\pm i/\varepsilon$, and so are 
$\alpha(P(z))$ and $\beta(P(z))$ in   view of (\ref{alpha-beta}).

Since by inductive hypothesis a.s. $\mbox{det}(K_{n-1})\neq 0$, by $\delta(P(z))\equiv 0$ we would get for MQ
$$
[-\mbox{det}(K_{n-1})(\Delta\phi_n(P(z)))^2+\beta(P(z))](1+(\varepsilon z)^2)^3
$$
$$
=-\varepsilon^4\mbox{det}(K_{n-1})(2+(\varepsilon z)^2)^2+\beta(P(z))(1+(\varepsilon z)^2)^3
$$
$$
\equiv - \alpha(P(z))\Delta\phi_n(P(z))(1+(\varepsilon z)^2)^3
$$
\begin{equation} \label{deltaMQ}
= - \varepsilon^2\alpha(P(z))(2+(\varepsilon z)^2)(1+(\varepsilon z)^2)^{3/2}\;,
\end{equation}
and for IMQ
$$
[-\mbox{det}(K_{n-1})(\Delta\phi_n(P(z)))^2+\beta(P(z))](1+(\varepsilon z)^2)^5
$$
$$
=\varepsilon^4\mbox{det}(K_{n-1})(2-(\varepsilon z)^2)^2+\beta(P(z))(1+(\varepsilon z))^5
$$
$$
\equiv -\alpha(P(z))\Delta\phi_n(P(z))(1+(\varepsilon z)^2)^5
$$
\begin{equation} \label{deltaIMQ}
=+\varepsilon^2\alpha(P(z))(2-(\varepsilon z)^2)(1+(\varepsilon z)^2)^{5/2}\;,
\end{equation}
which in both cases give a contradiction, because the first term 
in (\ref{deltaMQ}) and  (\ref{deltaIMQ}) is analytic at $\pm i/\varepsilon$, whereas the last has a branching point there. 

Then, 
$\mbox{det}(K_{n+1})=\delta(P_{n+1})$ is a.s. nonzero, 
by the already quoted fundamental result that the zero set of a not identically zero real analytic function on an open connected set $\Omega\subset \mathbb{R}^d$ is a null set for the Lebesgue measure (and thus also for any probability measure with density $\sigma\in L^1_+(\Omega)$). Indeed, 
denoting by $Z_\delta$ the zero set of $\delta$ in $\Omega$ and recalling that $\mbox{det}(K_{n-1})\neq 0$ (which a.s. holds) implies $\delta\not\equiv 0$, taking the probability of the corresponding events we get
$$ 
\mbox{prob}\{\mbox{det}(K_{n+1})=0\}=\mbox{prob}\{\delta(P_{n+1})=0\}
$$
$$
=\mbox{prob}\{\delta\equiv 0\}
+\mbox{prob}\{\delta\not\equiv 0\;\&\;P_{n+1}\in Z_\delta\}
=0+0=0\;,
$$ 
and the inductive step is completed.
\hspace{0.2cm} $\square$

\vskip0.5cm 
\noindent
{\bf Acknowledgements.} 

Work partially supported by the DOR funds of the University of Padova, and by the INdAM-GNCS 2024 Project “Kernel and polynomial methods for approximation and integration: theory and application software''. 

This research has been accomplished within the RITA ``Research ITalian network on Approximation", the SIMAI Activity Group ANA\&A, and the UMI Group TAA ``Approximation Theory and Applications".

\end{document}